\documentclass[12pt]{article}
\usepackage{hyperref}
\usepackage{import}
\usepackage{amsmath,amsthm,amsfonts}
\usepackage{amstext}
\usepackage{graphicx}
\usepackage{xcolor}
\usepackage{algorithmicx}
\usepackage[linesnumbered,ruled,vlined]{algorithm2e}
\usepackage{algpseudocode}
\usepackage[numbers]{natbib}
\usepackage[utf8]{inputenc} 
\usepackage[T1]{fontenc}    

\usepackage[top=2cm, bottom=2cm, left=2cm, right=2cm]{geometry}

\begin{document}

%
%
%

\newcommand{\nobj}{m}
\newcommand{\pfront}{\mathcal{P}}
\newcommand{\Prob}{\mathbb{P}}
\newcommand{\Yr}{\boldsymbol{\mathcal{Y}}}
\def\ds{\displaystyle}
\def\R{\mathbb{R}}
\newcommand{\Y}{\mathbf{Y}}
\newcommand{\F}{\mathbf{F}}
\newcommand{\z}{\mathbf{z}}
\newcommand{\s}{\mathbf{x}}
\newcommand{\Rset}{\mathbb{R}}
\newcommand{\yr}{\mathcal{Y}}
\newcommand{\fr}{\mathcal{F}}
\newcommand{\Fr}{\boldsymbol{\mathcal{F}}}
\newcommand{\xr}{\mathcal{X}}
\newcommand{\sr}{\mathbb{X}}
\newcommand{\esp}{\mathbb{E}}
\def\cov{\operatorname{cov}}

\newcommand{\x}{\mathbf{x}}
\newcommand{\y}{\mathbf{y}}
\newcommand{\uto}{\mathbf{u}}
\newcommand{\dis}{\mathbf{d}}
\newcommand{\Xset}{\mathcal{X}}
\newcommand{\PF}{\mathcal{P}}
\newcommand{\f}{\ensuremath{\mathbf{f}}}  
\newcommand{\M}{\ensuremath{\mathcal{M}}} 

\theoremstyle{definition}
\newtheorem{definition}{Definition}[section]
\newtheorem*{remark}{Remark}

\title{A game theoretic perspective on Bayesian multi-objective optimization}

\author{Micka\"{e}l Binois\thanks{Inria, Universit\'e C\^ote d'Azur, CNRS, LJAD, UMR 7351, Sophia Antipolis, France, mickael.binois@inria.fr, habbal@unice.fr}
    \and
  Abderrahmane Habbal\footnotemark[1]
    \and
  Victor Picheny\thanks{Secondmind, Cambridge, United Kingdom, victor@secondmind.ai}}

\date{}

\maketitle

\abstract{This chapter addresses the question of how to efficiently solve
many-objective optimization problems in a computationally demanding black-box
simulation context. We shall motivate the question by applications in machine
learning and engineering, and discuss specific harsh challenges in using
classical Pareto approaches when the number of objectives is four or more.
Then, we review solutions combining approaches from Bayesian optimization,
e.g., with Gaussian processes, and concepts from game theory like Nash
equilibria, Kalai-Smorodinsky solutions and detail extensions like
Nash-Kalai-Smorodinsky solutions.  We finally introduce the corresponding
algorithms and provide some illustrating results.}

\section{Introduction}

In the previous Chapter 11, methods for multi-objective (MO) Bayesian
optimization (BO) have been introduced. Scaling with many-objectives burst the
seams of these existing solutions on several aspects:
i) The positive-definiteness requirement for multi-output Gaussian process (GP) covariance matrix severely restricts their modelling ability, letting independent GPs as the default alternative that does not directly exploit possible correlations between objectives (see also Chapter 10);
ii) Summing the objectives is a go-to strategy often observed in practice. But picking weights in scalarization methods becomes more arbitrary, with effects that are difficult to apprehend. More generally, specifying preferences for a large number of  conflicting  objectives  can  lead  to  situations  where  the  connection  between  achievable solutions and the given preferences is not easy to interpret;
iii) For Pareto-based infill criteria, the computational cost grows quickly, for instance with hypervolume computations. Plus there are less closed form expressions available than for two or three objectives;
iv) The effect of more objectives on the acquisition function landscape is to be studied, but they are expected to be more multimodal due to larger Pareto sets.
And these are only a few additional technical problems faced by many-objective optimization (MaOO).

But even the desired outcome is problematic. As the number of objectives
increases, the dimension of the Pareto front generally grows accordingly, and
representing this latter is more complex. With a limited budget of
evaluations, this may not even be achievable: a few hundred design points
would not give an accurate representation of a complex ten-dimensional
manifold. In addition, even if representing the entire Pareto front was
achievable, the decision maker would be left with a very large set of
incomparable solutions, with limited practical relevance.

Hence, MaOO is more of an elicitation problem, for which good and principled
ways to choose solutions are needed. When the decision maker is involved in
the optimization process by giving feedback on pairs of solutions, methods
learning this preference information such as
\cite{Lepird2015,astudillo2020multi} are available. More diffuse preference
information can also be included via the choice of the reference point for
hypervolume computation (nullifying the contributions of equivalent or
dominated solutions), interactively modifying ranges for objectives
\cite{hakanen2017using} or via an aspiration point serving as target (see,
e.g., \cite{gaudrie2020targeting}). Ordering the importance of objectives is
another option, exploited e.g., in \cite{Abdolshah2019}. With many objectives,
possibly of various natures and scales, the burden on the decision maker is
heavier, and the effect on the outcome of these preferences is less intuitive.

When the decision is only made \textit{a posteriori}, knee points, Pareto
solutions from which improving slightly one objective would deteriorate much
more some others, are appealing to decision makers, or in the design of
evolutionary algorithms, see e.g., \cite{zhang2014knee}. But their existence
and location on the Pareto front may vary, especially in the MaOO case. Other
arguments have been directed toward finding a single solution at the
``center'' of the Pareto front to keep balance between objectives, e.g., by
\cite{gaudrie2018budgeted,paria2020flexible}. The rationale is that trade-off
solutions at the extremes are not desirable unless prior information has been
given. Finding such balanced solutions, without input from the decision maker
during the optimization, is the default framework considered in this chapter.
The difficulty lies in an appropriate and principled definition of a central
solution. A possible approach is to frame the MaOO problem as a game theory
one, where each objective is impersonated by a \textit{player}, and the
players bargain until they agree on a satisfying solution for all, referred to
as \textit{equilibrium}.

Game theory is already popular in a variety of applications such as
engineering \cite{desideri2012cooperation}, control \cite{habbal2013neumann},
multi-agent systems \cite{brown2015hierarchical} or machine learning
\cite{lanctot2012no}. Game theory has also revealed to be a powerful tool in
multi-disciplinary optimization, where decentralized strategies for
non-cooperative agents takes its full and relevant meaning, as referring to
discipline-specific solvers in a parallel, asynchronous and heterogeneous
multi-platform setting, see e.g., the pioneering framework developed in
\cite{periaux2015evolutionary}. It comes equipped with principled ways of
eliciting solutions (the \textit{bargaining axioms} of scale invariance,
symmetry, efficiency, and independence of irrelevant alternatives), echoing
several desirable properties for improvement functions used in MO BO
\cite{Svenson2011}: reflecting Pareto dominance, no input in the form of
external parameters (such as a reference point), or robustness to scale
invariance. Indeed, with more objectives of different nature, finding a
reasonable common scale for all is much harder and the effects cannot be
apprehended.

Compared to the more standard multi-objective BO framework described in the
previous chapter, the games framework comes with some specific challenges. The
definition of solutions such as Nash equilibrium or Kalai-Smorodinsky
solutions can be quite involved, and finding such equilibria generally implies
solving several inference tasks at once. For instance, this makes the use of
improvement-based approaches difficult. Still, several approaches have emerged
using regret or stepwise uncertainty reduction (SUR) frameworks, see e.g.,
\cite{picheny2019bayesian,Abdolshah2019} and references therein.

The structure of the chapter is as follows. In Section \ref{sec:game}, we
briefly review key concepts in game theory such as Nash equilibria. Then in
Section \ref{sec:algo} we present high-level BO approaches to solve
many-objective problems under the Games paradigm. Illustrations are provided
in Section \ref{sec:applications} before discussing remaining challenges and
perspectives in Section \ref{sec:perspectives}.

\section{Game equilibria to solution elicitation}
\label{sec:game}

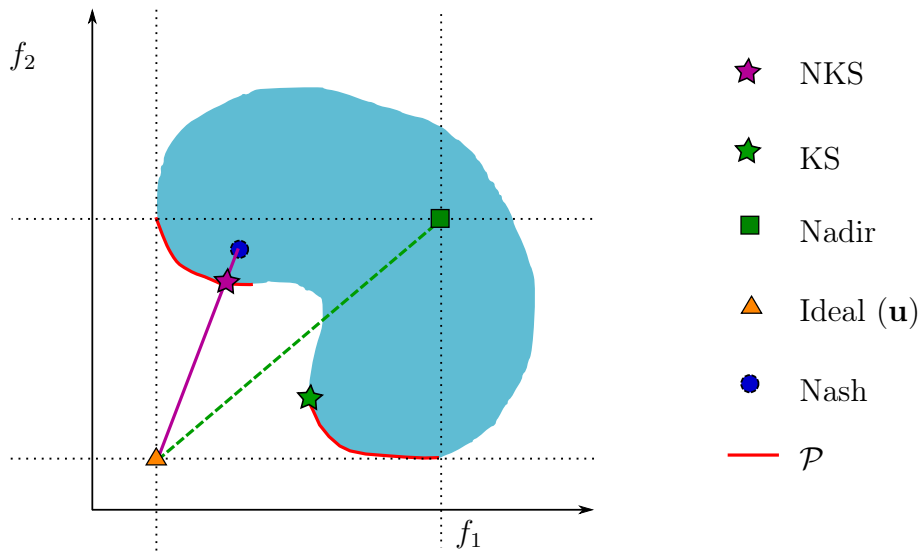
\begin{figure}[htpb]%
  \centering
  \def\svgwidth{0.85\textwidth}
\begingroup%
  \makeatletter%
  \providecommand\color[2][]{%
    \errmessage{(Inkscape) Color is used for the text in Inkscape, but the package 'color.sty' is not loaded}%
    \renewcommand\color[2][]{}%
  }%
  \providecommand\transparent[1]{%
    \errmessage{(Inkscape) Transparency is used (non-zero) for the text in Inkscape, but the package 'transparent.sty' is not loaded}%
    \renewcommand\transparent[1]{}%
  }%
  \providecommand\rotatebox[2]{#2}%
  \newcommand*\fsize{\dimexpr\f@size pt\relax}%
  \newcommand*\lineheight[1]{\fontsize{\fsize}{#1\fsize}\selectfont}%
  \ifx\svgwidth\undefined%
    \setlength{\unitlength}{265.8013916bp}%
    \ifx\svgscale\undefined%
      \relax%
    \else%
      \setlength{\unitlength}{\unitlength * \real{\svgscale}}%
    \fi%
  \else%
    \setlength{\unitlength}{\svgwidth}%
  \fi%
  \global\let\svgwidth\undefined%
  \global\let\svgscale\undefined%
  \makeatother%
  \begin{picture}(1,0.48635228)%
    \lineheight{1}%
    \setlength\tabcolsep{0pt}%
    \put(0,0){\includegraphics[width=\unitlength,page=1]{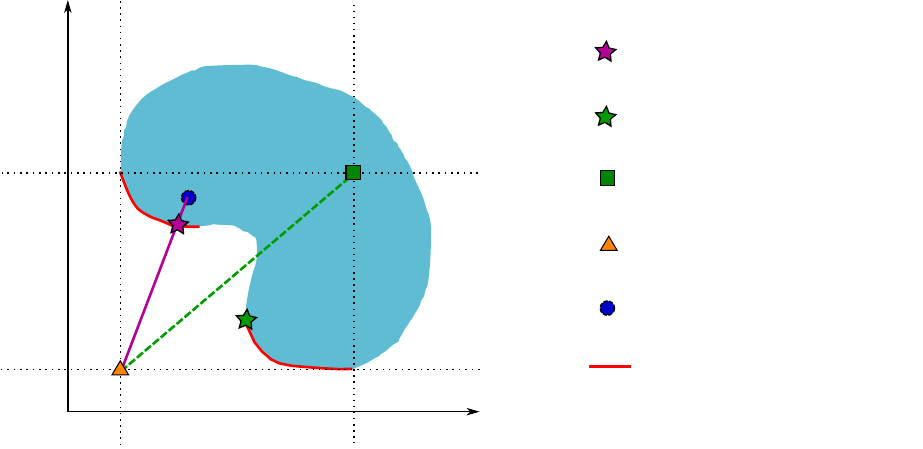}}%
    \put(-0.00072322,0.43681358){\color[rgb]{0,0,0}\makebox(0,0)[lt]{\lineheight{0}\smash{\begin{tabular}[t]{l}$f_2$ \end{tabular}}}}%
    \put(0.39614194,0.01162466){\color[rgb]{0,0,0}\makebox(0,0)[lt]{\lineheight{0}\smash{\begin{tabular}[t]{l}$f_1$ \end{tabular}}}}%
    \put(0.70222166,0.41873614){\color[rgb]{0,0,0}\makebox(0,0)[lt]{\lineheight{0}\smash{\begin{tabular}[t]{l}NKS \end{tabular}}}}%
    \put(0.70222166,0.20800631){\color[rgb]{0,0,0}\makebox(0,0)[lt]{\lineheight{0}\smash{\begin{tabular}[t]{l}Ideal ($\mathbf{u})$ \end{tabular}}}}%
    \put(0.70222166,0.34324583){\color[rgb]{0,0,0}\makebox(0,0)[lt]{\lineheight{0}\smash{\begin{tabular}[t]{l}KS \end{tabular}}}}%
    \put(0.70222166,0.27834775){\color[rgb]{0,0,0}\makebox(0,0)[lt]{\lineheight{0}\smash{\begin{tabular}[t]{l}Nadir \end{tabular}}}}%
    \put(0.70222166,0.13748645){\color[rgb]{0,0,0}\makebox(0,0)[lt]{\lineheight{0}\smash{\begin{tabular}[t]{l}Nash \end{tabular}}}}%
    \put(0.70222166,0.07661195){\color[rgb]{0,0,0}\makebox(0,0)[lt]{\lineheight{0}\smash{\begin{tabular}[t]{l}$\mathcal{P}$ \end{tabular}}}}%
  \end{picture}%
\endgroup%
\caption{Illustration of solutions of two different elicitation problems with
two different disagreement points $\dis$: either the Nadir point or a Nash
equilibrium, resulting in the KS and NKS points (stars) on the Pareto front
$\pfront$, respectively. The shaded area shows the feasible objective space.}
\label{fig:disagrement}
\end{figure}

The standard multi-objective optimization problem (MOP) corresponds to the
simultaneous optimization of all objectives:
\begin{equation}
    (MOP) \quad \min \limits_{\x \in \Xset} (f_1(\x), \dots, f_m(\x)).
    \label{eq:MOP}
\end{equation}
Besides considering objectives as individual players' goals, general
non-cooperative games, and related solutions, such as Nash ones, do need
territory splitting. That is, partitioning of the optimization variables (the
input space) among players. Other games and related solution concepts depend
on choices of anchor points, as is the case for the Kalai-Smorodinsky (KS)
solution, which depends on ideal and disagreement points, as illustrated in
Figure \ref{fig:disagrement}. We shall introduce in the following the Nash
equilibrium concept, then motivated by the generic inefficiency of the latter,
move to considering the Kalai-Smorodinsky solution.

\subsection{Nash games and equilibria}

When considering primarily the standard (static, under complete information)
Nash equilibrium problem \cite{Gibbons-Game-92}, each objective becomes a
player's outcome. Compared to a standard MOP (\ref{eq:MOP}) where players
share the control of the same set of variables (or \emph{action space}),
variables are uniquely allocated to a player in the so-called territory
splitting. Denote $\x_i$ the variables of player $i$ and $\xr_i$ its
corresponding action space, where $\xr = \prod_i \xr_i$. Accordingly, the
variable vector $\x$ consists of block components $\x_1, \ldots, \x_m$
$\left(\x = (\x_j)_{1\le j \le m}\right)$. We shall use the convention
$f_i(\x) = f_i(\x_i, \x_{-i})$ when we need to emphasize the role of $\x_i$,
where $\x_{-i}$ is the subset of variables controlled by players $j \neq i$.

\begin{definition}
A Nash equilibrium problem (NEP) consists of $m \ge 2$ decision makers (i.e.,
players), where each player $i \in \{1, \ldots, m\}$ wants to solve its
optimization problem:
\begin{equation}\label{eq:NEP}
 (\mathcal{P}_i) \quad \min_{\x_i \in \xr_i} f_i(\x),
\end{equation}
and $\mathbf{f}(\x) = \left[f_1(\x),\ldots, f_m(\x)\right]: \xr \subset \Rset^n
\rightarrow \Rset^m$ (with $n \geq m$) denotes a vector of cost functions
(a.k.a. pay-off or utility functions when maximized), $f_i$ denotes the specific cost
function of player $i$.
\end{definition}

\begin{definition}
 A Nash equilibrium $\x^* \in \xr$ is a strategy such that:
\begin{equation} \label{eq:NEdef}
(NE) \quad \forall i,\ 1\le i \le m, \quad  \x^*_i  = \arg \min_{\x_i \in \xr_i } f_i(\x_i, \x^*_{-i}).
\end{equation}
\end{definition}
In other words, when all players have chosen to play a NE, then no single
player has incentive to move from his $\x^*_i$.\\

The above definition shows some attractive features of NE: it is
scale-invariant, does not depend on arbitrary parameters (such as the
reference point), and convey a notion of first-order stationarity. On the
other hand, it requires partitioning the inputs between players, which may
come naturally for some problems, but may be problematic for many MaOO
problems. Plus an arbitrary partitioning can have an unknown outcome on the
solution. One option proposed in \cite{desideri2014} is to define the
partitioning based on sensitivity analysis of one main objective. Hence a
perspective would be to define an optimal partitioning in the sense that it
allocates variables (or a linear combination of variables) to the player it
has most influence on.

Besides the need to define a territory splitting, it is important to notice
that, generically, Nash equilibria are not efficient, i.e., do not belong to
the underlying set of best compromise solutions, the so-called Pareto front,
of the objective vector $(f_i(\x))_{\x \in \Xset}$. Indeed, NE efficiency may
happen when one of the players has control on all the optimization variables
while the others have control on nothing (a somehow degenerate territory
splitting). Still, NE have a couple of advantages: a nice notion of
stationarity, they are generally well balanced and, of interest for many
objectives, scale invariant.
\begin{remark}
When block components $\x_1, \ldots, \x_m$ of the variable vector $\x$ overlap
due to, e.g., constraints, then one should modify (\ref{eq:NEP}) and
(\ref{eq:NEdef}) accordingly, so that $\x_i \in \xr_i (\x_{-i}) $. In the
literature, the Nash equilibrium problem is then referred to as a generalized
Nash equilibrium problem (GNEP).
\end{remark} 

\begin{remark}
The NE solution is scale invariant, and more generally invariant under any
strictly increasing transformation $\Psi_i$ ($ 1\le i \le m$) since, in this
case, (\ref{eq:NEdef}) is equivalent to
\begin{equation*}
 \forall i,\ 1\le i \le m, \quad  \Psi_i(f_i(\x^*_i, \x^*_{-i})) \le \Psi_i (f_i(\x_i, \x^*_{-i}))
\end{equation*}

\end{remark}

Interestingly, it could be shown at least in the convex case that if players
share probabilistic control on the optimization variable components (i.e., any
player $j$  may control the same component $\x_i$ with probability $p_{ij}$)
then there exists a probability matrix $(p_{ij})$ such that the associated
Nash equilibrium (with an {\it ad hoc} definition) lies in the Pareto front,
see \cite{aboulaich:hal-00648693} for a sketch of the approach.

From a practical viewpoint, computing NE for games stated in continuous
variable settings (as opposed to discrete or finite games, {e.g.,} in vector
spaces with Banach or Hilbert structure) can be based on variational analysis,
e.g., with the classical fixed-point algorithms to solve NEPs
\cite{uryas1994relaxation,MR942837,MR899829}. A modified notion of
Karush-Kuhn-Tucker (KKT) points, adapted to (generalized) NE, is developed by
\cite{kanzow2016augmented}, to propose a dedicated augmented Lagrangian
method. Yet these methods require too many evaluations to tackle directly
expensive black-boxes, requiring specific BO algorithms that are presented in
Section \ref{sec:algo}.

\subsection{The Kalai-Smorodinsky solution}\label{sec:KS}

The Kalai-Smorodinsky solution was first proposed by Kalai and Smorodinsky in
1975 as an alternative to the Nash bargaining solution in cooperative
bargaining. Differently from Nash equilibria, no partitioning of the decision
variable among players is required by the KS solution concept.

The problem is as follows: starting from a {\em disagreement} or {\em status
quo} point $\dis$ in the objective space, the players aim at maximizing their
own benefit while moving from $\dis$ toward the Pareto front (i.e., the
efficiency set). The KS solution is of egalitarian inspiration
\cite{conley1991bargaining} and states that the selected efficient solution
should yield equal benefit ratio to all the players. Indeed, given the utopia
(or ideal, or shadow) point $\mathbf{u}\in\mathbb{R}^\nobj$ defined by
\[u_{i} = \min_{\x \in \xr}\ f_{i}(\x),
\]
selecting any compromise solution $\mathbf{y} = [f_{1}(\x), \ldots,
f_{m}(\x)]$ would yield, for objective $i$, a benefit ratio
\[r_{i}(\x)= \frac{d_{i} - y_{i}}{d_{i} - u_{i}}.\]

Notice that the benefit from staying at $\dis$ is zero, while it is maximal
for the generically unfeasible point $\mathbf{u}$. The KS solution is the
Pareto optimal choice $\mathbf{y}^* = [f_{1}(\x^*), \ldots, f_{m}(\x^*)]$ for
which all the benefit ratios $r_{i}(\x)$ are equal:
\begin{equation*}
    \x^* \in \left\{\x \in \Xset \text{~s.t.~} r_{1}(\x) = \dots = r_{m}(\x) \right\} \cap \left\{\x \in \Xset \text{~s.t.~} \mathbf{f}(\x) \in \pfront \right\}.
\end{equation*}
Geometrically, $\mathbf{y}^*$ is the intersection point of the Pareto front
and the line $(\dis,\uto)$ (see Figure \ref{fig:disagrement}), which may be
empty. For $m=2$, $\mathbf{y}^*$ exists if the Pareto front is continuous. For
the general case, we use here the extension of the KS solution proposed by
\cite{hougaard2003nonconvex} under the name \textit{efficient maxmin
solution}. Since the intersection with the $(\dis,\uto)$ line might not be
feasible, there is a necessary trade-off between Pareto optimality and
centrality. The efficient maxmin solution is defined as the Pareto-optimal
solution that maximizes the smallest benefit ratio among players, that is:
\begin{equation*}\label{eq:defKS}
  \mathbf{x}^{**} \in  \arg \max_{\mathbf{y} \in \pfront} \min_{1 \leq i \leq \nobj} r_{i}(\x).
\end{equation*}
It is straightforward that when the intersection is feasible, then
$\mathbf{y}^{*}$ and $\mathbf{y}^{**} := \mathbf{f}(\x^{**})$ coincide. Figure
\ref{fig:disagrement} shows $\mathbf{y}^{**}$ in the situation when the
feasible space is nonconvex. $\mathbf{y}^{**}$ is always on the Pareto front
(hence not necessarily on the $(\dis,\uto)$ line). In the following, we refer
indifferently to $\mathbf{y}^{*}$ (if it exists) and $\mathbf{y}^{**}$ as the
KS solution. Note that this definition also extends to discrete Pareto sets.

For $\nobj = 2$, the (non-extended) KS solution can be axiomatically
characterized as the unique solution fulfilling all the bargaining solution
axioms, which are: Pareto optimality, symmetry, affine invariance, and
restricted monotonicity \cite{KSE1975}. For $\nobj \ge 3$, there is no
axiomatic setting, the KS solution being required to fulfill Pareto
optimality, affine invariance, and equity in benefit ratio. Still, KS is
particularly attractive in a many-objective context since it scales naturally
to a large number of objectives and returns a single solution, avoiding the
difficulty of exploring and approximating large $\nobj$-dimensional Pareto
fronts--especially with a limited number of observations.

Closely related to the Kalai and Smorodinsky solution is the reference point
methodology, developed by \cite{wierzbicki1980use}. The reference point
methodology uses achievement functions, which refer to aspiration and
reservation references. Ideal and nadir points are particular instances of
such reference points, respectively, and achievement functions play the same
role as the benefit ratio of our present setting. Moreover, the notion of
neutral compromise and max-min (of achievement functions) were introduced by
Wierzbicki in the cited reference, yielding a framework objectively close to
the one of Kalai and Smorodinsky. One important difference between the
reference point and KS approaches is that the latter relies on a game
theoretic axiomatic construction, in view of palliating the Nash equilibrium
inefficiency.

\subsection{Disagreement point choice}

Clearly, the KS solution strongly depends on the choice of the disagreement
point $\dis$. A standard choice is the nadir point $N$ given by $N_i =\max_{\x
\in \text{Pareto set}}\ f_{(i)}(\x)$. Some authors introduced alternative
definitions, called extended KS, to alleviate the critical dependence on the
choice of $\dis$ \cite{bozbay2012bargaining}, for instance by taking as
disagreement point the Nash equilibrium arising from a previously played
non-cooperative game. But such a choice would need a pre-bargaining split of the
decision variable vector $\x$ among the $\nobj$ players. When some relevant territory
splitting is agreed on, then Nash yields KS an interesting and non arbitrary
disagreement point, and in return, KS yields the Nash game an interesting and
non arbitrary efficient solution, starting from the non efficient equilibrium.
The resulting Nash-Kalai-Smorodinsky solution is denoted NKS.

The nadir point is very useful to rescale the Pareto front in case objectives
are of different nature or scale, see e.g., \cite{bechikh2010estimating}. As
such, it makes a natural disagreement solution. Still, finding the nadir point
is a complex task, especially under a limited budget, as it involves exploring
regions good for a few objectives with worse values on others. The many
objective context makes this issue more prominent. It is related to the
problem of Pareto resistance, i.e., regions of the space with good values on a
few objectives but not Pareto optimal (see, e.g., \cite{fieldsend2016enabling}
and references therein). Unless some objectives are more important than
others, the corresponding extremal regions of the Pareto front are of little
interest for selecting a single good solution.

A simpler disagreement point is via the \emph{pseudo-nadir} point, defined as
the worst objective solution over designs achieving the minimum for one
objective: $\tilde{N}_i = \max_{\x \in \{ \x^{(1)*}, \dots, \x^{(m)*} \} }\
f_{i}(\x)$, where $\x^{(j)*} \in \arg \min_{\x \in \Xset} f_j(\x)$. It can be
extracted directly from the pay-off table, but it is not robust as it can
under- or over-estimate the true nadir point, see e.g.,
\cite{miettinen2012nonlinear} for details. While for two objectives, nadir and
pseudo-nadir coincide, this is not the case anymore in higher dimensions.
Using this pseudo nadir point removes the need to search for extremal
non-dominated solutions on the Pareto front, and aligns with finding the ideal
point.

The egalitarian inspiration of the KS solution is based on the assumption that
all the objectives have equal importance. However, in many practical cases,
the end user may want to favor a subset of objectives which are of primary
importance, while still retaining the others for optimization. One way of
incorporating those \textit{preferences}
\cite{thiele2009preference,wang2017mini,purshouse2014review} is to discard
solutions with extreme values and actually solve the following problem
(constrained version of the original MOP (\ref{eq:MOP})):
\begin{eqnarray*}
 \min_{\x \in \Xset} & \left\{ f_{1}(\x), \ldots, f_{\nobj}(\x)
 \right\}\nonumber \\ \text{s.t.} & f_{i}(\x) \leq c_i, \qquad i \in J
 \subset \left[1, \ldots, \nobj \right],\label{eq:constraints}
\end{eqnarray*}
with $c_i$'s predefined (or interactively defined) constants. 
Choosing a tight
value (i.e., difficult to attain) for $c_i$ may discard a large portion of the
Pareto set and favor the i-th objective over the others.

Incorporating such preferences in the KS solution can simply be done by using
$\mathbf{c}$ as the disagreement point if all objectives are constrained, or
by replacing the coordinates of the nadir with the $c_i$ values. In a game
theory context, this would mean that each player would state a limit of
acceptance for his objective before starting the cooperative bargaining. From
geometrical considerations, one may note that $\dis$ (and hence $\mathbf{c}$)
does not need to be a feasible point. As long as $\dis$ is dominated by the
utopia point $\uto$, the KS solution remains at the intersection of the Pareto
front and the $(\dis,\uto)$ line. An example where preferences from a decision
maker are integrated this way is provided in \cite{binois2019kalai}.

As shown for instance in \cite{binois2019kalai}, by applying a monotonic
rescaling to the objectives, the corresponding KS solution can move along the
Pareto front. An example of a monotonic value function applied to a raw
objective corresponds to taking the logarithm of one objective for a physical
or modeling reason. Since Nash equilibrium are invariant to such rescaling,
they are more stable than nadir or pseudo-nadir based KS solutions, but still
affected. One proposed option to alleviate this effect is to use ranks instead
of objective values directly. Interestingly, it induces a dependency on the
parameterization of the problem. Studying these points can include the use of
a measure on the input space and the use of copulas in the objective one,
\cite{binois2019kalai,Binois2015b}, with further work needed to quantify the
effect on the Pareto front. The combination with more involved preference
learning methods such as \cite{Lepird2015,astudillo2020multi} could also be
investigated.

\section{Bayesian optimization algorithms for games}
\label{sec:algo}

Solving expensive black-box games is an emerging topic in Bayesian
optimization
\cite{aprem2018bayesian,al2018approximating,picheny2019bayesian,binois2019kalai}.
Similarly to the single and multi-objective cases, BO leverages the
probabilistic information available from the Gaussian process models to
balance exploration and exploitation in searching for the game equilibrium.
Other metamodelling options could be entertained, as long as they also provide
mean and variance predictions at any location, plus the ability to get
realizations (posterior draws) at arbitrary locations. We present here several
algorithms to solve such equilibrium or solution finding problems.

For the Nash equilibrium (Section \ref{sec:fixedpoint}),
\cite{aprem2018bayesian} proposed an upper confidence bound-like acquisition
function based on a regret function approximating the game-theoretic regret.
The case of potential games \cite{gonzalez2016survey}, when a single function
can summarize the game, commonly arising in problems with shared resources, is
considered by \cite{al2018approximating}, with a dedicated expected
improvement-like acquisition function. In \cite{picheny2019bayesian,
binois2019kalai}, the flexible and general step-wise uncertainty reduction is
applied to estimate different classes of equilibria, including Nash and KS,
which we present in Section \ref{sec:sur}. Although not proposed in the
literature (up to our knowledge), a related family of acquisition functions
based on Thompson sampling is also applicable, as we describe briefly in
Section \ref{sec:ts}.

\subsection{Fixed point approaches for the Nash equilibrium}\label{sec:fixedpoint}

Fixed points methods, see an example in Algorithm \ref{alg:fixed-point}, are
one way to get Nash equilibrium. In the expensive case, a na\"ive version is
simply to replace the true objectives by their GP mean prediction:
$\hat{f}_i$, $1 \leq i \leq m$. While computationally efficient and possibly
sufficient for a coarse estimation, such an approach does not address the
exploration-exploitation trade-off, and may lead to poor estimates if the
models are not accurate representations of the objectives.

\begin{algorithm}[!ht]
\caption{Pseudo-code for the fixed-point approach \cite{uryas1994relaxation}}
{\bf Require: } $m$: number of players,  $0<\alpha<1$ : relaxation factor, $n_{\max}$ : max iterations\;
Construct initial strategy $\x^{(0)}$\;
\While{$n \leq n_{\max}$}{
  Compute in parallel: $\forall i,\ 1\le i \le m, \quad \z^{(k+1)}_{i} = \arg \min_{\x_i \in \Xset_i} \hat{f}_i(\x^{(k)}_{-i} ,\x_{i})$\;
  Update : $\x^{(k+1)} = \alpha \z^{(k+1)} + (1-\alpha) \x^{(k)} $\;
  \If{$\| \x^{(k+1)} - \x^{(k)}\|$ small enough}{
      exit\;
  }  
} 
{\bf Ensure: } For all $i=1,...,m$, $\x^{*}_{i} = \arg \min_{\x_i \in \Xset_i} J_i(\x^{*}_{-i}, \x_{i})$
\label{alg:fixed-point}
\end{algorithm}

\cite{al2018approximating} provides an acquisition function based on the upper
confidence bound to estimate Nash equilibrium for continuous games. They focus
on the game theoretic regret, that is the most any player can gain from
deviating from the current solution. Following small deviations for each
player in turns is also the underlying principle in fixed point methods. The
corresponding acquisition function aims to select as next design to evaluate
one achieving the minimum of the regret of the GP surrogate. But rather than
optimizing directly to get the minimum objective values, these are replaced by
an approximation. That is, for a player $i$ at $\x_i$, the minimum value is
replaced by the sum of the mean objective value over $\x_{-i}$ with the
corresponding standard deviation, scaled by a factor $\gamma$. These can be
computed analytically for GPs for some separable covariance functions.

In the special case of potential games \cite{gonzalez2016survey}, also related
to fixed point methods, the strategy of improving as much as possible one
player leads to the Nash equilibrium (and not any fixed point). This is
exploited by \cite{aprem2018bayesian} who directly considers the potential. By
definition, the potential $\Phi:\Xset \to \R$ is a function such that
$f_i(\x', \x_{-i}) - f_i(\x'', \x_{-i}) = \Phi(\x',
\x_{-i}) - \Phi(\x'', \x_{-i})$, $\forall i, \x', \x''$. As the available
measurements are on the $f_i$'s and not the potential directly, constructing the GP of the potential $\Phi$ involves taking into account integral and gradient operators. Once the GP is fitted, they proposed an acquisition function similar to expected improvement \cite{Mockus74}.

\subsection{Stepwise uncertainty reduction}\label{sec:sur}

Step-wise uncertainty reduction (SUR) methods, see, e.g.,
\cite{geman1996active,Villemonteix2009,chevalier2014fast}, require an
uncertainty measure $\Gamma$ about a quantity of interest, here the game
equilibrium. Then new evaluations are selected sequentially in order to reduce
this uncertainty, leading to an accurate identification of the quantity of
interest. This framework is well-adapted to games as it allows to integrate
various learning tasks into a single goal (for instance, the definition of KS
type solutions involves several other unknown quantities, such as the ideal
and nadir points).

For any multivariate function $\mathbf{f}$ defined over $\Xset$, with image
$\mathbb{Y}$, denote by $\Psi: \mathbb{Y} \rightarrow \Rset^\nobj$ the mapping
that associates its equilibrium (Nash, KS or other). When the multivariate
function is replaced by a GP emulator conditioned on $n$ observations $\left(
\mathbf{f}(\x^{(1)}), \dots, \mathbf{f}(\x^{(n)}) \right)$, $\Psi$ is defined
on the corresponding distribution $\Y_n()$. Then $\Psi(\Y_n)$ is a random
vector (of unknown distribution). Loosely speaking, the spread of the
corresponding distribution in the objective space reflects the uncertainty on
the solution location, given by the uncertainty measure $\Gamma(\Y_n)$. One
measure of variability of a vector is the determinant of its covariance matrix
\cite{fedorov1972theory}: $\Gamma(\Y_n) = \det \left[\cov \left(\Psi(\Y_n)
\right) \right]$, while other information theoretic measures could be used.
For instance, this could be the conditional entropy as in
\cite{Villemonteix2009}.

The SUR strategy greedily chooses the next observation that reduces the most
this uncertainty:
\begin{equation*}
 \max_{\x \in \Xset} \left\{ \Gamma(\Y_n) - \Gamma(\Y_{n, \x}) \right\},
\end{equation*}
where $\Y_{n, \x}$ is the GP conditioned on $\{ \mathbf{f}(\x^{(1)}), \ldots,  \mathbf{f}(\x^{(n)}),
\mathbf{f}(\x)\}$. However, such an ideal strategy is not tractable as it would
require evaluating all $\mathbf{f}(\x)$ while maximizing over $\Xset$.

A more manageable strategy is to consider the expected uncertainty reduction:
$\Gamma(\Y_n) - \esp_{\Y_n(\x)} \left[ \Gamma(\Y_{n, \x}) \right]$, where
$\esp_{\Y_n(\x)}$ denotes the expectation taken over $\Y_n(\x)$. Removing the
constant term $\Gamma(\Y_n)$, the policy is defined with
\begin{equation*}\label{eq:nextpoint}
  \x \in \arg \min_{\x \in \Xset} \left\{ J(\x) = \esp_{\Y_n(\x)} \left[ \Gamma(\Y_{n, \x}) \right] \right\}.
\end{equation*}
Informally speaking, the variation of $\Psi$ is partly caused by not knowing
the precise localization of the solution of interest, plus eventually from not
knowing the values of the anchor points. The realizations corresponding to
these specific points would either be for designs close to the solution of the
GP means, or for designs with large variance. Hence $J(\x)$ defines a
trade-off between exploration and exploitation, as well as a trade-off between
the different learning tasks ($\dis$ and $\uto$ points, Nash equilibrium,
Pareto front). Hence the BO optimization loop is fully defined (see Chapter
11, Algorithm 1) and the difficulty boils down to evaluating $J$ efficiently.
 
A crucial aspect in enabling the use of the SUR strategy is the ability to
generate realizations or estimate $\Y_n$ efficiently. When $\Xset$ is discrete
(or can be discretized), the approaches proposed in
\cite{picheny2019bayesian,binois2019kalai} relies on the use of conditional
simulations, i.e., joint posterior samples on $\Xset$ or a well chosen subset
of it, coupled with fast update formulas for the resulting ensembles,
following \cite{chevalier2015fast}. This approach hinges on the discretization
size used, which needs to remain in the thousands due to the cubic
computational complexity in the number of samples. The targeted solution is
then computed on each realization, resulting in samples from $\Y_n$. Note that
for discrete Nash equilibrium computations, all combinations of the selected
strategies for each player must be simulated at. An alternative solution,
applicable to continuous $\Xset$ would rely on approximated sample paths with
a closed form expression as used, e.g., by \cite{Mutny2018}. Again, the
appropriate solution can be obtained from the samples, which is inexpensive
compared to using the expensive black-box.

The resulting criterion for KSE estimation automatically interleaves improving
the estimation of the specific quantities involved, by balancing between
estimation of the Nash equilibrium, of the objective-wise minima for the ideal
point and of the intersection with the Pareto front. A different use of the
posterior distribution is also possible, with Thompson sampling.

\subsection{Thompson sampling}\label{sec:ts}

Over the past decade, Thompson sampling (TS) \cite{Thompson1933} has become a
very popular algorithm to solve sequential optimization problems. In a
nutshell, TS proceeds by sequentially sampling from the posterior distribution
of the problem solution, allowing for efficiently addressing
exploration-exploitation trade-offs. In a GP-based setting, this is simply
achieved by sampling from the posterior distribution of the objectives and
choosing as the next sampling point the input that realizes the equilibrium on
the sample.

As it would not directly try to pinpoint the ideal point (say), it cannot be
directly transposed for estimating KS solutions. Still, it may serve to
estimate Nash equilibria. However, up to our knowledge this solution has not
been studied yet (the closest proposition is an algorithm in
\cite{picheny2019bayesian} based on the probability of realizing the optimum).
Note that the same approximation of the samples as for SUR can be used. Next
we illustrate the results of this latter on a practical application.

\section{Application example: engineering test case}
\label{sec:applications}

We consider the switching ripple suppressor design problem for voltage source
inversion in powered systems, which was proposed by \cite{zhang2019switching},
and also tested in \cite{he2019repository}. The device is composed of three
components, one with three inductors, one with parallel LC resonant branches
and a capacitor branch. Arguing that this type of problem is naturally many
objective, \cite{zhang2019switching} describes several types of conflicting
objectives: harmonic attenuations of the switching ripple suppressor under
different frequencies, power factor, resonance suppression and inductor cost.

In the proposed problem, where $n_r$ is the number of resonant branches, there
are $n_r+4$ variables and $n_r+1$ objectives: variables are inductors $L_1,
L_2, L_3$ that can be related to the inductor cost objective ($f_{n_r + 1}$),
and variables $C_1, \dots, C_{n_r}$ related to the goal of attenuating the
harmonics at resonant frequencies ($f_1, \dots, f_{n_r}$). The last variable
is $C_f$. Note that five additional design constraints can also be added. The
cost of evaluation is here negligible, but it showcases a realistic
application example where the true Pareto front is unknown.

\begin{figure}[htbp]
\centering
 \includegraphics[trim=0mm 0mm 0mm 0mm, width=\textwidth, clip]{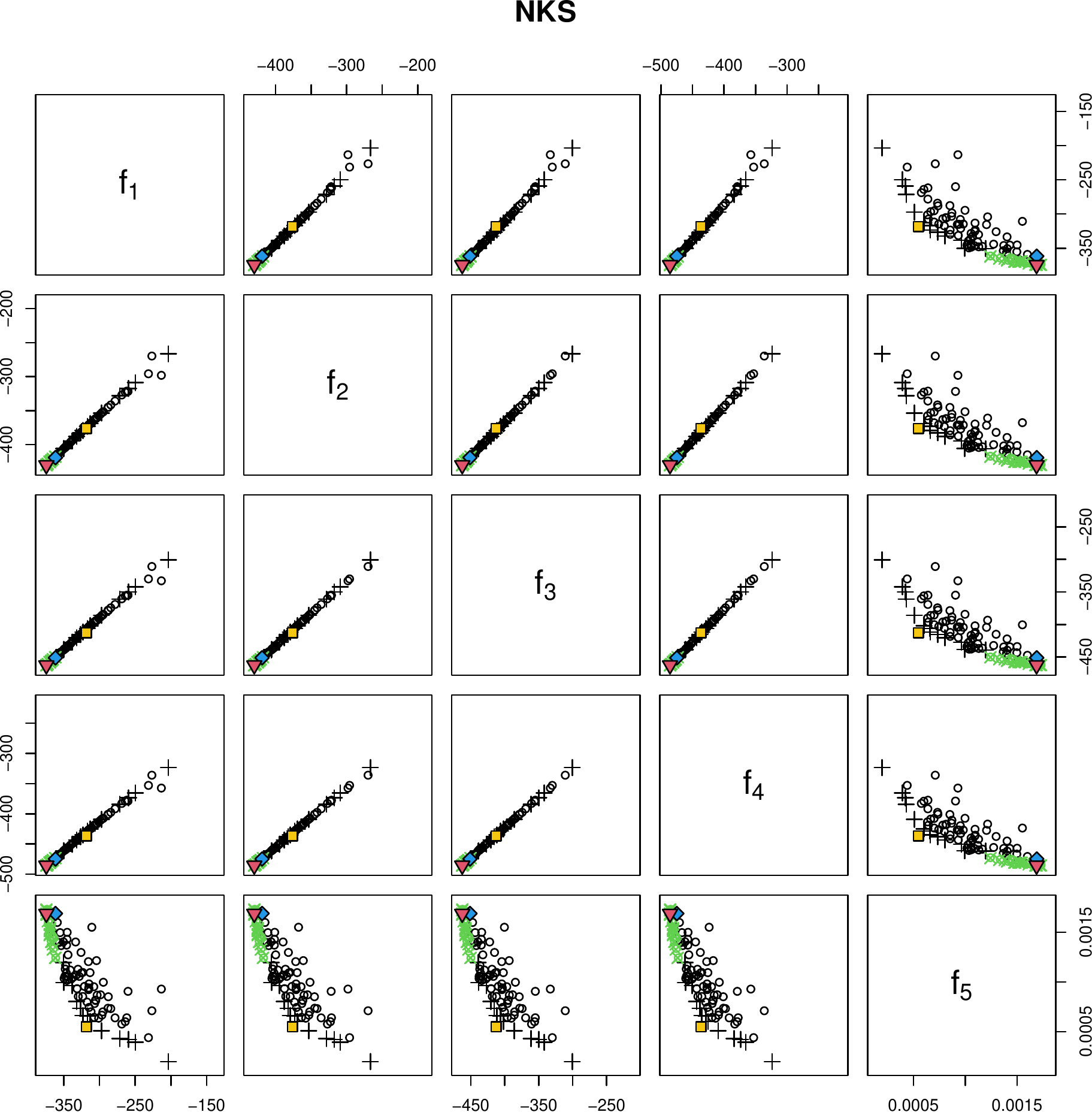}
 \caption{Scatter plot matrix of the objectives obtained when estimating the NKS solution by the SUR strategy. 
 Black circles: dominated initial designs, black +: nondominated initial
 designs, green crosses: sequentially added solutions (circles if dominated),
 blue triangles: estimated Nash equilibrium, and red diamond: estimated NKS
 solution. The KS solution (yellow square) is marked for reference only.}
\label{fig:rippleNKS}
\end{figure}

\begin{figure}[htbp]
\centering
 \includegraphics[trim=20mm 15mm 15mm 15mm, width=\textwidth, clip]{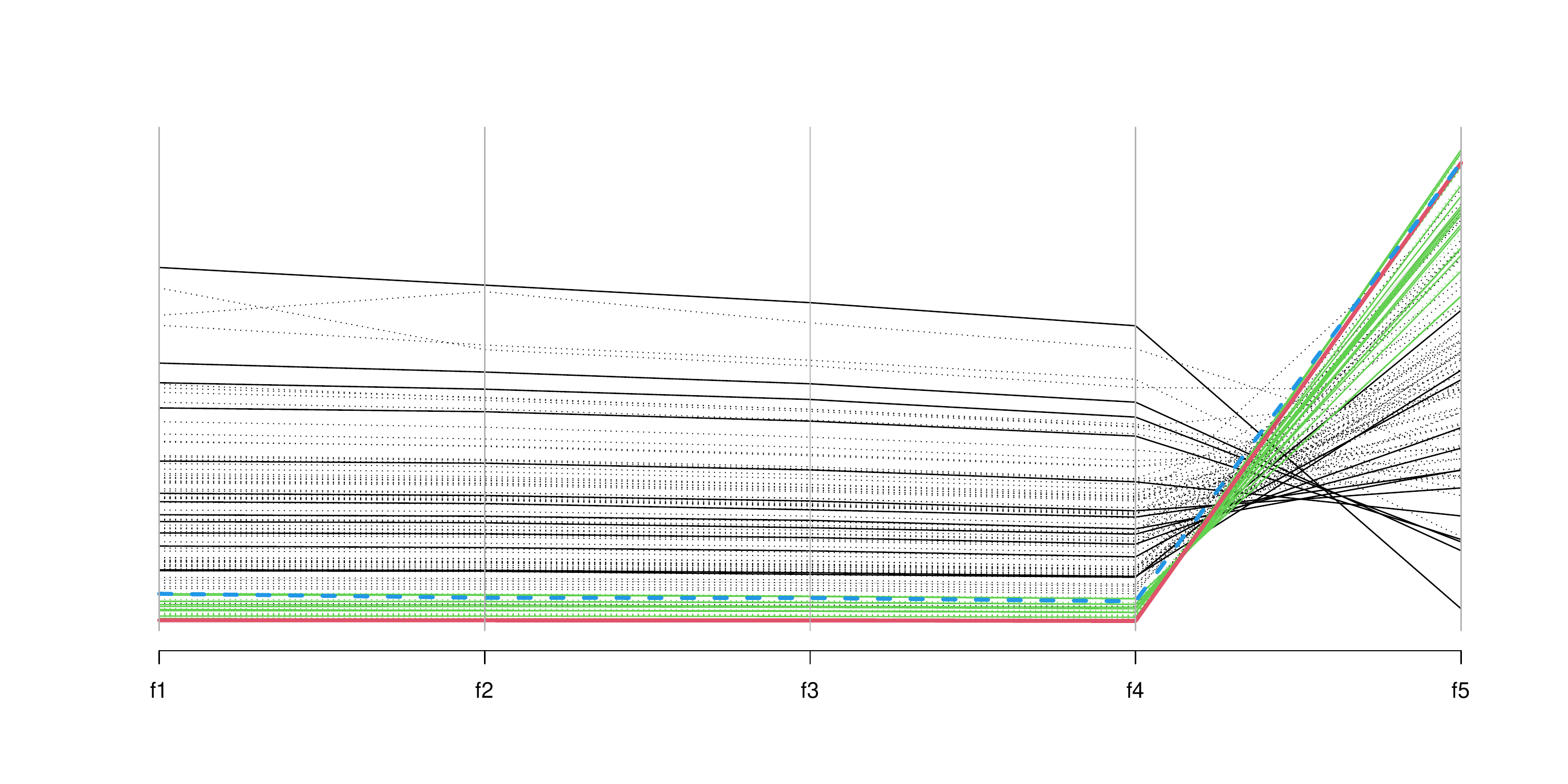}
 \caption{Parallel coordinates plot corresponding to the objectives of Figure \ref{fig:rippleNKS}. 
 Colors are the same, dominated solutions are in dotted lines while the Nash
 equilibrium is marked by the large blue dashed line and the estimated NKS
 solution by the large red line.}
\label{fig:parrippleNKS}
\end{figure}

We start by estimating the Nash-Kalai-Smorodinsky equilibrium in the $n_r = 4$
case; with the following partitioning: $(f_1 : C_1, C_f), (f_2 : C_2), (f_3 :
C_3), (f_4: C_4), (f_5 : L_1, L_2, L_3)$. For the SUR strategy, we discretize
each input subspace using Latin hypercube samples of sizes (26, 11, 11, 11,
51) and use a budget of 80 evaluations for the initial design followed by 70
sequential iterations. Note that this discretization of the input space with
the partitioning results in $\approx 2\times10^6$ possible combinations for
the Nash equilibrium. To keep it tractable for the SUR procedure, filtering is
performed to select 1000 solutions to actually draw posterior realizations,
before evaluating the criterion on the 200 most promising candidates. More
details on the procedure can be found in
\cite{picheny2019bayesian,binois2019kalai}. The result obtained using the
\texttt{GPGame} package \cite{Picheny2020} is described in Figures
\ref{fig:rippleNKS} and \ref{fig:parrippleNKS}. It appears that the first four
objectives are highly correlated, which does not affect the estimation of the
solution. Most of the initial designs (67/80) are dominated, while this is the
case for two thirds (47/70) of the added ones, possibly due to the search of
the inefficient Nash equilibrium. These new designs are all concentrated in
the same region of the output space. The NE is relatively close to the Pareto
front in this case, with the obtained solution dominating it.

\begin{figure}[htbp]
\centering
 \includegraphics[trim=0mm 0mm 0mm 0mm, width=\textwidth, clip]{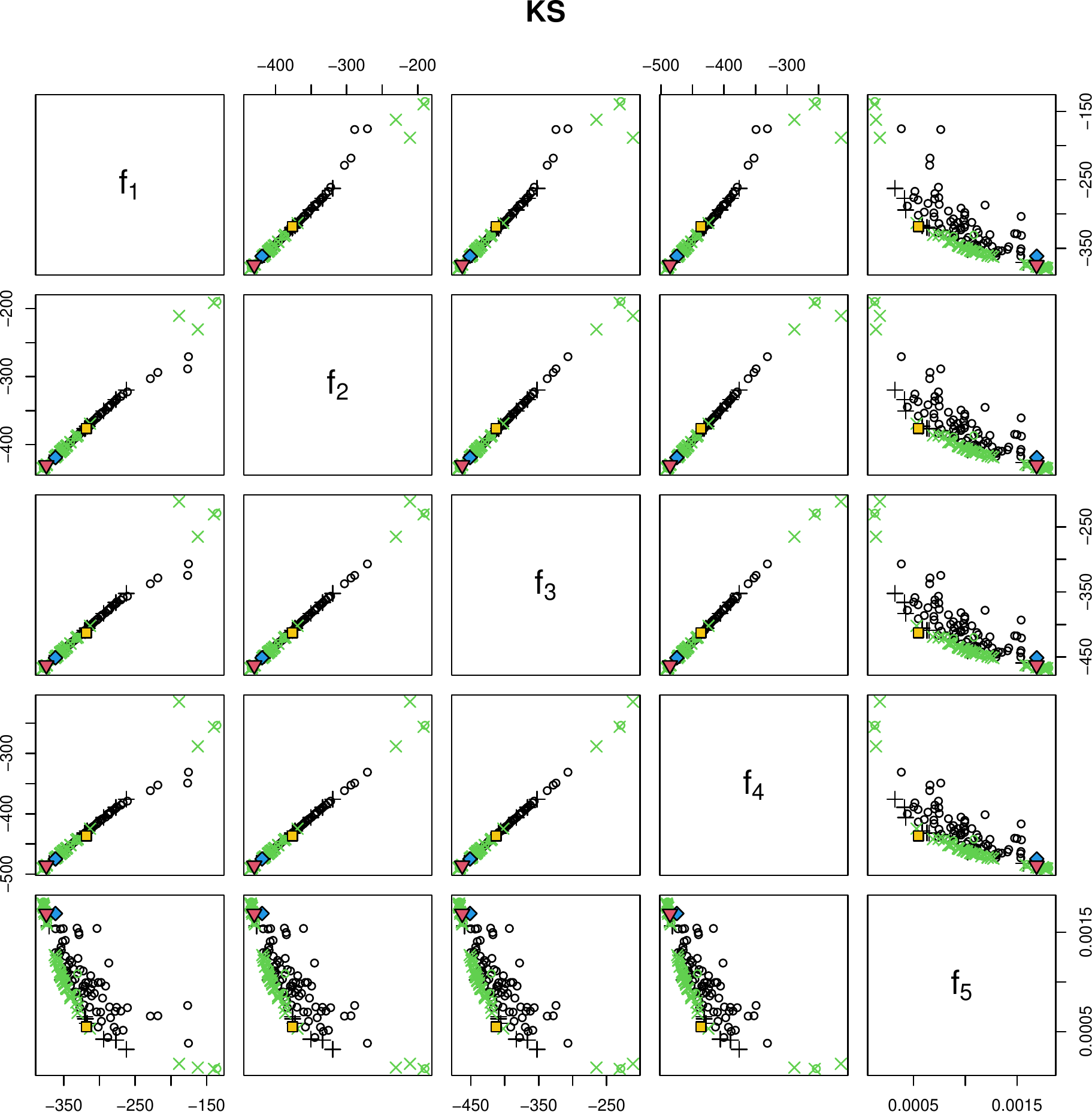}
 \caption{Scatter plot matrix of the objectives obtained when estimating the KS solution by the SUR strategy. The legend is as for Figure
 \ref{fig:rippleNKS} except that the yellow square indicates the KS solution.
 Here the Nash solution (blue diamond) and NKS solution (red triangle) are
 marked for reference only.}
\label{fig:rippleKS}
\end{figure}

\begin{figure}[htbp]
\centering
 \includegraphics[trim=20mm 15mm 15mm 15mm, width=\textwidth, clip]{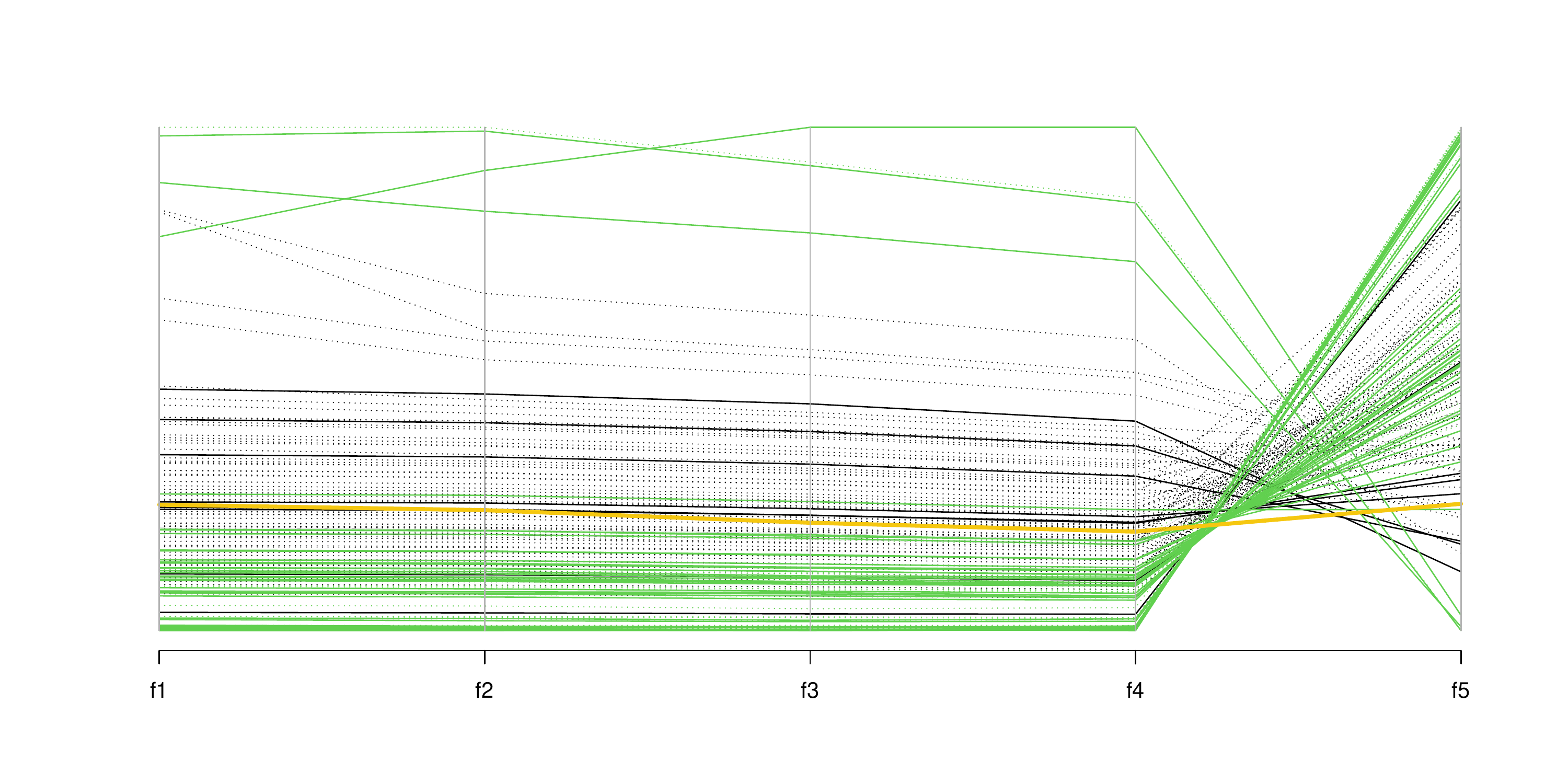}
 \caption{Parallel coordinates plot corresponding to the objectives of Figure
 \ref{fig:rippleKS}. Colors are the same, dominated solution are depicted with
 dotted lines while the estimated KS solution is marked by the large yellow
 line.}
\label{fig:parrippleKS}
\end{figure}

Then we consider the regular Kalai-Smorodinsky solution, whose results are in
Figures \ref{fig:rippleKS} and \ref{fig:parrippleKS}. This time the input
space is discretized with $10^6$ possible solutions, uniformly sampled. Again,
most initial design are dominated (72/80), but this is the case for only a
third (28/70) of the sequentially added ones. Figure \ref{fig:parrippleKS}
shows that the KS solution is almost a straight line, an indication of the
centrality of this solution. It can be seen that, in the sequential procedure,
new points are added at the extremities of the Pareto front to locate the
nadir, as well as in the central part.

Choosing between the obtained NKS and KS solutions now depends on the
viewpoint. The spread of the latter is much broader in the objective space
than for NKS. Note that the same range is used on all the Figures. To further
ease comparison, we added the KS solution on Figure \ref{fig:rippleNKS} and
the Nash and NKS ones on Figure \ref{fig:rippleKS}, but there were not part of
the solutions evaluated in their corresponding BO loops, nor their targets. In
case the range of objectives is relevant to the decision maker, the centrality
and equity of benefits on all objectives may be appealing to the decision
maker. If instead the territory splitting for the objectives is important, the
relative scaling of the objectives irrelevant, then the NKS solution would be
preferred (also over staying on the Nash equilibrium). As a perspective, it
would be pertinent to additionally take into account the constraints defined
in \cite{zhang2019switching}, for which the Pareto front becomes more complex.

\section{What is done and what remains}
\label{sec:perspectives}

We reviewed the enrichment of many-objective optimization with concepts from
game theory. In the context of limited budgets of evaluations when the
extensive search of the entire Pareto front is unreachable, game theory can
guide the estimation of a single solution chosen for its centrality and
possible additional game theoretic properties. Several BO-based methodologies
have been developed for estimating the corresponding solution, which can scale
to many-objective setups. As this is an emerging topic in BO, it seems that
many algorithmic developments could be proposed, for instance based on
Thompson sampling or on entropy search.

For the Nash and NKS solutions, a partitioning of the input space is
necessary, while it is seldom included in the multi-objective context. A
perspective is to randomize this partitioning when it does not come directly
from the problem definition. The possible advantage would be to focus on a set
of different Nash equilibria, all with stationarity properties, providing more
than a single solution.

As with the engineering example we presented, there are usually additional
constraints in the problem definition. In the SUR procedure described above,
dealing with constraints could amount to filtering unfeasible solutions on the
realizations. Other options could include the use of the probability of
feasibility as in \cite{Schonlau1998}. For more details on the handling of
constraints in BO, we refer to \cite{Letham2018} and references therein. Then
the difference between an objective and a constraint is sometimes fuzzy.
Especially when the constraints can be modified, taking them as objectives
instead, akin to multiobjectivization \cite{knowles2001reducing}, could help
finding feasible solutions (or work in the case the constrained problem is
unfeasible). Conversely, prioritizing some objectives
\cite{desideri:hal-02285197} while defining constraints on others can help
control the trade-offs when moving on the Pareto front.


\bibliographystyle{plain}


\end{document}